\documentclass[11pt]{amsart}
\usepackage{amssymb}
\usepackage{amsmath}
\usepackage[active]{srcltx}
\usepackage{t1enc}
\usepackage[latin2]{inputenc}
\usepackage{verbatim}
\usepackage{amsmath,amsfonts,amssymb,amsthm}
\usepackage[mathcal]{eucal}
\usepackage{enumerate}
\usepackage[centertags]{amsmath}
\usepackage{graphics}

\newtheorem{theorem}{Theorem}

\newtheorem{lemma}{Lemma}

\begin{document}
\author{George Tephnadze}
\title[Fejér means]{Approximation by Walsh-Kaczmarz-Fejér means on the Hardy
space }
\address{G. Tephnadze, Department of Mathematics, Faculty of Exact and
Natural Sciences, Tbilisi State University, Chavchavadze str. 1, Tbilisi
0128, Georgia and Department of Engineering Sciences and Mathematics, Lule%
\aa {} University of Technology, SE-971 87, Lule\aa {}, Sweden.}
\email{giorgitephnadze@gmail.com}
\thanks{The research was supported by Shota Rustaveli National Science
Foundation grant no.13/06 (Geometry of function spaces, interpolation and
embedding theorems).}
\date{}
\maketitle

\begin{abstract}
The main aim of this paper is to find necessary and sufficient conditions
for the convergence of Walsh-Kaczmarz-Fejér means in the terms of the
modulus of continuity on the Hardy spaces $H_{p},$ when $0<p\leq 1/2.$
\end{abstract}

\date{}

\section{Introduction}

\textbf{2000 Mathematics Subject Classification.} 42C10.

\textbf{Key words and phrases:} Walsh-Kaczmarz system, Fejér means,
martingale Hardy space, modulus of continuity.

The a.e. convergence of \ Walsh-Fejér means $\sigma _{n}f$ was proved by
Fine \cite{F}. Yano \cite{Yano} has proved that $\left\Vert K_{n}\right\Vert
\leq 2$ ($n\in \mathbb{N}$).Consequently, $\left\Vert \sigma
_{n}f-f\right\Vert _{p}\rightarrow 0$, as $n\rightarrow \infty $ ($f\in
L_{p},$ $1\leq p\leq \infty $). However (see \cite{JOO, SWSP}) the rate of
convergence can not be better then $O\left( n^{-1}\right) $ $\left(
n\rightarrow \infty \right) $ for nonconstant functions. Fridli \cite{FR}
use dyadic modulus of continuity to characterize the set of functions in the
space $L_{p}$, whose Walsh-Fejér means converge at the given rate.

In 1975 Schipp \cite{Sch} showed that the maximal operator $\sigma ^{\ast }$
is of weak type $(1,1)$ and of type $(p,p)$ for $1<p\leq \infty $. The
boundedness fails to hold for $p=1$. But, Fujii \cite{Fu} proved that $%
\sigma ^{\ast }$ is bounded from the dyadic Hardy space $H_{1}$ to the space
$L_{1}$. The theorem of Fujii was extended by Weisz \cite{We4}, he showed
that the maximal operator $\sigma ^{\ast }$ is bounded from the martingale
Hardy space $H_{p}$ to the space $L_{p}$ for $p>1/2.$ Simon gave a
counterexample \cite{Si}, which shows that the boundedness does not hold for
$0<p<1/2.$ The counterexample for $p=1/2$ due to Goginava \cite{GogiEJA}
(see also \cite{tep1}). In the endpoint case $p=1/2$ Weisz \cite{We5} proved
that $\sigma ^{\ast }$ is bounded from the Hardy space $H_{1/2}$ to the
space weak-$L_{1/2}$.

In 1948 $\breve{\text{S}}$neider \cite{Snei} introduced the Walsh-Kaczmarz
system and showed that the inequality
\begin{equation*}
\limsup_{n\rightarrow \infty }\frac{D_{n}^{\kappa }(x)}{\log n}\geq C>0
\end{equation*}%
holds a.e. In 1974 Schipp \cite{Sch2} and Young \cite{Y} proved that the
Walsh-Kaczmarz system is a convergence system. In 1981, Skvortsov \cite{Sk1}
showed that the Fejér means with respect to the Walsh-Kaczmarz system
converge uniformly to $f$ for any continuous functions $f$. Gát \cite{gat}
proved that, for any integrable functions, the Fejér means with respect to
the Walsh-Kaczmarz system converge almost everywhere to the function. He
showed that the maximal operator of Walsh-Kaczmarz-Fejér means $\sigma
^{\kappa ,\ast }$ is weak type $(1,1)$ and of type $(p,p)$ for all $1<p\leq
\infty $. Gát's result was generalized by Simon \cite{S2}, who showed that
the maximal operator $\sigma ^{\kappa ,\ast }$ is of type $(H_{p},L_{p})$
for $p>1/2$. In the endpoint case $p=1/2$ Goginava \cite{Gog-PM} (see also
\cite{GNCz,GGN-SSMH}) proved that maximal operator $\sigma ^{\kappa ,\ast }$
is not of type $(H_{1/2},L_{1/2})$ and Weisz \cite{We5} showed that the
maximal operator is of weak type $(H_{1/2},L_{1/2})$. In \cite{GNCz} and
\cite{tep4} it is proved that the maximal operator $\widetilde{\sigma }$ $%
^{\kappa ,\ast \,}$defined by
\begin{equation*}
\widetilde{\sigma }_{p}^{\kappa ,\ast }:=\sup_{n\in \mathbb{N}}\frac{%
\left\vert \sigma _{n}^{\kappa }\right\vert }{\left( n+1\right) ^{1/p-2}\log
^{2\left[ 1/2+p\right] }\left( n+1\right) },
\end{equation*}%
where $0<p\leq 1/2$ and $\left[ 1/2+p\right] $ denotes integer part of $%
1/2+p,$ is bounded from the Hardy space $H_{p}$ to the space $L_{p}$. As a
corollary we get
\begin{equation}
\left\Vert \sigma _{n}^{\kappa }f\right\Vert _{p}\leq c_{p}\left( n+1\right)
^{1/p-2}\log ^{2\left[ 1/2+p\right] }\left( n+1\right) \left\Vert
f\right\Vert _{H_{p}}.  \label{c1}
\end{equation}

Moreover, for any nondecreasing function $\varphi :\mathbb{N}_{+}\rightarrow
\lbrack 1,$ $\infty )$ satisfying the condition
\begin{equation}
\overline{\lim_{n\rightarrow \infty }}\frac{\left( n+1\right) ^{1/p-2}\log
^{2\left[ 1/2+p\right] }\left( n+1\right) }{\varphi \left( n\right) }%
=+\infty ,  \label{cond1}
\end{equation}%
there exists a martingale $f\in H_{p},$ such that
\begin{equation*}
\underset{n\in \mathbb{N}}{\sup }\left\Vert \frac{\sigma _{n}^{\kappa }f}{%
\varphi \left( n\right) }\right\Vert _{p}=\infty .
\end{equation*}

For Walsh system analogical theorems are proved in \cite{GoSzeged}, \cite%
{tep2} and \cite{tep3}.

It is well-known (see \cite{we2}) that
\begin{equation}
\Vert \sigma _{2^{n}}^{\kappa }(f)-f\Vert _{H_{p}}\rightarrow 0,\text{ \
when }n\rightarrow \infty ,\text{ \ }\left( f\in H_{p},\text{ }p>0\right) .
\label{31}
\end{equation}

Móricz and Siddiqi \cite{Mor} investigates the approximation properties of
some special Nörlund means of Walsh-Fourier series of $L_{p}$ function in
norm. Fridli, Manchanda and Siddiqi \cite{fms} improved and extended the
results of Móricz and Siddiqi \cite{Mor} among them in $H_{p}$ norm, where $%
0<p<1.$ In \cite{gog8} Goginava investigated the behavior of Cesáro means of
Walsh-Fourier series in detail. Analogical results for
Walsh-Kaczmarz-Fourier series were proved by Nagy in \cite{N1,N2,N3,N4}.

The main aim of this paper is to find necessary and sufficient conditions
for the convergence of Fejér means in the terms of the modulus of continuity
on the Hardy spaces $H_{p},$ when $0<p<1/2$ and $p=1/2,$ respectively. That
is, the most investigated end point case $p=1/2$ is also discussed.

The paper is organized as following: In Section 3 we present and discuss the
main results and in Section 4 the auxiliary propositions will be given. In
Section 5 the the proofs can be found. Moreover, in order not to disturb our
discussions in these Sections some definitions and notations are given in
Section 2.

\section{Definitions and Notations}

Now, we give a brief introduction to the theory of dyadic analysis \cite%
{SWSP}. Let $\mathbb{N}_{+}$ denote the set of the positive integers, $%
\mathbb{N}:=\mathbb{N}_{+}\cup \{0\}.$ Denote ${\mathbb{Z}}_{2}$ the
discrete cyclic group of order 2, that is ${\mathbb{Z}}_{2}=\{0,1\},$ where
the group operation is the modulo 2 addition and every subset is open. The
Haar measure on ${\mathbb{Z}}_{2}$ is given such that the measure of a
singleton is 1/2. Let $G$ be the complete direct product of the countable
infinite copies of the compact groups ${\mathbb{Z}}_{2}.$ The elements of $G$
are of the form $x=\left( x_{0},x_{1},...,x_{k},...\right) $ with $x_{k}\in
\{0,1\}\left( k\in \mathbb{N}\right) .$ The group operation on $G$ is the
coordinatewise addition, the measure (denoted by $\mu $) and the topology
are the product measure and topology. The compact Abelian group $G$ is
called the Walsh group. A base for the neighborhoods of $G$ can be given in
the following way:%
\begin{equation*}
I_{0}\left( x\right) :=G,
\end{equation*}%
\begin{equation*}
I_{n}\left( x\right) :=I_{n}\left( x_{0},...,x_{n-1}\right) :=\left\{ y\in
G:\,y=\left( x_{0},...,x_{n-1},y_{n},y_{n+1},...\right) \right\} ,
\end{equation*}%
where $x\in G$ and $n\in \mathbb{N}.$ These sets are called dyadic
intervals. Let $0=\left( 0:i\in \mathbb{N}\right) \in G$ denote the null
element of $G,$ and $I_{n}:=I_{n}\left( 0\right) \,\left( n\in \mathbb{N}%
\right) .$ Set $e_{n}:=\left( 0,...,0,x_{n}=1,0,...\right) \in G$ and%
\begin{equation*}
J_{N}^{m,l}:=I_{N}\left( x_{0},...,x_{m}=1,0,...,x_{l}=1,0,...,0\right) ,%
\text{ \ }l=0,...,N-1,\text{ \ }m=-1,...,l.
\end{equation*}

For $k\in \mathbb{N}$ and $x\in G$ denote
\begin{equation*}
r_{k}\left( x\right) :=\left( -1\right) ^{x_{k}}
\end{equation*}%
the $k$-th Rademacher function.

If $n\in \mathbb{N}$, then $n=\sum\limits_{i=0}^{\infty }n_{i}2^{i}$ can be
written, where $n_{i}\in \{0,1\}\quad \left( i\in \mathbb{N}\right) $, i. e.
$n$ is expressed in the number system of base 2. Denote $\left\vert
n\right\vert :=\max \{j\in \mathbb{N}\mathbf{:}n_{j}\neq 0\}$, that is $%
2^{\left\vert n\right\vert }\leq n<2^{\left\vert n\right\vert +1}.$

The Walsh-Paley system is defined as the sequence of Walsh-Paley functions:
\begin{equation*}
w_{n}\left( x\right) :=\prod\limits_{k=0}^{\infty }\left( r_{k}\left(
x\right) \right) ^{n_{k}}=r_{\left\vert n\right\vert }\left( x\right) \left(
-1\right) ^{\sum\limits_{k=0}^{\left\vert n\right\vert -1}n_{k}x_{k}}\quad
\left( x\in G,\text{ }n\in \mathbb{N}\right) .
\end{equation*}

The Walsh-Kaczmarz functions are defined by $\kappa _{0}=1$ and for $n\geq 1$

\begin{equation*}
\kappa _{n}\left( x\right) :=r_{\left| n\right| }\left( x\right)
\prod\limits_{k=0}^{\left| n\right| -1}\left( r_{\left| n\right| -1-k}\left(
x\right) \right) ^{n_{k}}=r_{\left| n\right| }\left( x\right) \left(
-1\right) ^{\sum\limits_{k=0}^{\left| n\right| -1}n_{k}x_{_{\left| n\right|
-1-k}}}.
\end{equation*}

Skvortsov (see \cite{Sk1}) gave a relation between the Walsh-Kaczmarz
functions and Walsh-Paley functions by the of the transformation $\tau
_{A}:G\rightarrow G$ defined by

\begin{equation*}
\tau _{A}\left( x\right) :=\left(
x_{A-1},x_{A-2},...,x_{1},x_{0,}x_{A},x_{A+1},...\right) ,
\end{equation*}%
for $A\in \mathbb{N}.$ By the definition we have

\begin{equation*}
\kappa _{n}\left( x\right) =r_{\left\vert n\right\vert }\left( x\right)
w_{n-2^{\left\vert n\right\vert }}\left( \tau _{\left\vert n\right\vert
}\left( x\right) \right) ,\qquad \left( n\in \mathbb{N}\mathbf{,}\text{ }%
x\in G\right) .
\end{equation*}

The Dirichlet kernels are defined

\begin{equation*}
D_{n}^{\alpha }:=\sum_{k=0}^{n-1}\alpha _{k\text{ }},
\end{equation*}%
where $\alpha _{n\text{ }}=w_{n}$ or $\kappa _{n}$ $\left( n\in \mathbb{N}%
\right) ,$ $D_{0}^{\alpha }:=0.$ the $2^{n}$-th Dirichlet kernels have a
closed form (see e.g. \cite{SWSP})

\begin{equation}
D_{2^{n}}^{w}(x)=D_{2^{n}}^{\kappa }(x)=D_{2^{n}}(x)=\left\{
\begin{array}{ll}
2^{n} & \text{\thinspace \thinspace }\,\text{if\thinspace \thinspace
\thinspace\ }x\in I_{n}, \\
0 & \text{ \thinspace if}\,\ \,x\notin I_{n}%
\end{array}%
\right.  \label{3}
\end{equation}%
The norm (or quasi-norm) of the space $L_{p}(G)$ is defined by \qquad

\begin{equation*}
\left\Vert f\right\Vert _{p}:=\left( \int_{G}\left\vert f\right\vert
^{p}d\mu \right) ^{1/p}\qquad \left( 0<p<\infty \right) .
\end{equation*}%
The space $L_{p,\infty }(G)$ consists of all measurable functions $f$ for
which
\begin{equation*}
\left\Vert f\right\Vert _{L_{p,\infty }}:=\underset{\lambda >0}{\sup }%
\lambda \mu \left\{ \left\vert f\right\vert >\lambda \right\} ^{1/p}\leq
c<\infty .
\end{equation*}

The $\sigma $-algebra generated by the dyadic intervals of measure $2^{-k}$
will be denoted by $F_{k}$ $\left( k\in \mathbb{N}\right) .$ Denote by $%
f=\left( f^{\left( n\right) },n\in \mathbb{N}\right) $ a martingale with
respect to $\left( F_{n},n\in \mathbb{N}\right) $ (for details see, e. g.
\cite{we2}). The maximal function of a martingale $f$ is defined by
\begin{equation*}
f^{\ast }=\sup\limits_{n\in \mathbb{N}}\left\vert f^{\left( n\right)
}\right\vert .
\end{equation*}

In case $f\in L_{1}\left( G\right) $, the maximal function can also be given
by
\begin{equation*}
f^{\ast }\left( x\right) =\sup\limits_{n\in \mathbb{N}}\frac{1}{\mu \left(
I_{n}(x)\right) }\left\vert \int\limits_{I_{n}(x)}f\left( u\right) d\mu
\left( u\right) \right\vert ,\ \ x\in G.
\end{equation*}

For $0<p<\infty $ the Hardy martingale space $H_{p}(G)$ consists of all
martingales for which

\begin{equation*}
\left\Vert f\right\Vert _{H_{p}}:=\left\Vert f^{\ast }\right\Vert
_{p}<\infty .
\end{equation*}

The best approximation of $f\in L_{p}(G)$ ($1\leq p\in \infty $) is defined
as%
\begin{equation*}
E_{n}\left( f,L_{p}\right) =\inf_{\psi \in \emph{p}_{n}}\left\Vert f-\psi
\right\Vert _{p}
\end{equation*}%
where $\emph{p}_{n}$ is set of all Walsh-Kaczmarz polynomials of order less
than $n\in \mathbb{N}$.

The integrated modulus of continuity of the function $f\in L_{p}(G)$, is
defined by

\begin{equation*}
\omega _{p}\left( 1/2^{n},f\right) =\sup\limits_{h\in I_{n}}\left\Vert
f\left( \cdot +h\right) -f\left( \cdot \right) \right\Vert _{p}.
\end{equation*}

The concept of modulus of continuity in $H_{p}(G)$ $\left( 0<p\leq 1\right) $
is defined in the following way
\begin{equation*}
\omega _{H_{p}}\left( 1/2^{n},f\right) :=\left\Vert f-S_{2^{n}}f\right\Vert
_{H_{p}}.
\end{equation*}

Since $\left\Vert f\right\Vert _{H_{p}}\sim \left\Vert f\right\Vert _{p}$,
when $p>1$, we obtain that
\begin{equation*}
\omega _{H_{p}}\left( 1/2^{n},f\right) \sim \left\Vert
f-S_{2^{n}}f\right\Vert _{p},\text{ \ \ }p>1.
\end{equation*}

On the ather hand Watari \cite{wat} showed that there are strong connection
among $\omega _{p}\left( 1/2^{n},f\right) ,$ $E_{2^{n}}\left( f,L_{p}\right)
$ and $\left\Vert f-S_{2^{n}}f\right\Vert _{p},$ \ \ $p\geq 1,$ $n\in
\mathbb{N}$. In particular%
\begin{equation*}
\omega _{p}\left( 1/2^{n},f\right) /2\leq \left\Vert f-S_{2^{n}}f\right\Vert
_{p}\leq \omega _{p}\left( 1/2^{n},f\right) ,
\end{equation*}%
and%
\begin{equation*}
\left\Vert f-S_{2^{n}}f\right\Vert _{p}/2\leq E_{2^{n}}\left( f,L_{p}\right)
\leq \left\Vert f-S_{2^{n}}f\right\Vert _{p}.
\end{equation*}

If $f\in L_{1}\left( G\right) $, then it is easy to show that the sequence $%
\left( S_{2^{n}}f:n\in \mathbb{N}\right) $ is a martingale. If $f$ is a
martingale, that is $f=(f^{\left( 0\right) },f^{\left( 1\right) },...)$ then
the Walsh-(Kaczmarz)-Fourier coefficients must be defined in a little bit
different way:
\begin{equation*}
\widehat{f}\left( i\right) =\lim\limits_{k\rightarrow \infty
}\int\limits_{G}f^{\left( k\right) }\left( x\right) \alpha _{i}\left(
x\right) d\mu \left( x\right) ,\ \ (\alpha _{i}=w_{i}\text{ or }\kappa _{i}).
\end{equation*}

The Walsh-(Kaczmarz)-Fourier coefficients of $f\in L_{1}\left( G\right) $
are the same as the ones of the martingale $\left( S_{2^{n}}f:n\in \mathbb{N}%
\right) $ obtained from $f$.

The partial sums and Fejér means of the Walsh-(Kaczmarz)-Fourier v of
martingale $f\in H_{p}$ are defined by
\begin{equation*}
S_{n}^{\alpha }(f):=\sum\limits_{i=0}^{n-1}\widehat{f}\left( i\right) \alpha
_{i},\quad \sigma _{n}^{\alpha }(f)=\frac{1}{n}\sum\limits_{j=1}^{n}S_{j}^{%
\alpha }(f),\text{ \ \ \ \ }(\alpha =w\text{ or }\kappa )
\end{equation*}%
respectively.

The Fejér kernel of order $n$ of the Walsh-(Kaczmarz)-Fourier series defined
by
\begin{equation*}
K_{n}^{\alpha }:=\frac{1}{n}\sum\limits_{k=1}^{n}D_{k}^{\alpha },\text{
\qquad }(\alpha =w\text{ or }\kappa ).
\end{equation*}

A bounded measurable function $a$ is p-atom, if there exists an interval $I$%
, such that

\begin{equation*}
\int_{I}ad\mu =0,\text{ \ \ \ \ \ }\left\Vert a\right\Vert _{\infty }\leq
\mu \left( I\right) ^{-1/p},\text{ \ \ \ \ supp}\left( a\right) \subset I.
\end{equation*}

For the martingale $f=\sum_{n=0}^{\infty }\left( f_{n}-f_{n-1}\right) $ the
conjugate transforms are defined as
\begin{equation*}
\widetilde{f^{\left( t\right) }}=\overset{\infty }{\underset{n=0}{\sum }}%
r_{n}\left( t\right) \left( f_{n}-f_{n-1}\right) ,
\end{equation*}%
where $t\in G$ is fixed. Note that $\widetilde{f^{\left( 0\right) }}=f.$ As
is well known (see \cite{we2})
\begin{equation}
\left\Vert \widetilde{f^{\left( t\right) }}\right\Vert _{H_{p}}=\left\Vert
f\right\Vert _{H_{p}},\text{ \ \ \ }\left\Vert f\right\Vert _{H_{p}}^{p}\sim
\int_{G}\left\Vert \widetilde{f^{\left( t\right) }}\right\Vert _{p}^{p}dt,%
\text{ \ \ \ \ \ \ }\widetilde{\sigma _{n}^{\kappa }f^{\left( t\right) }}%
=\sigma _{n}^{\kappa }\widetilde{f^{\left( t\right) }}.  \label{5.1}
\end{equation}

\section{Formulation of Main Results}

\begin{theorem}
a) Let $0<p<1/2$ and
\begin{equation*}
\omega _{H_{p}}\left( 1/2^{k},f\right) =o\left( 1/2^{k\left( 1/p-2\right)
}\right) ,\text{ as }k\rightarrow \infty .
\end{equation*}%
Then
\begin{equation}
\left\Vert \sigma _{n}^{\kappa }\left( f\right) -f\right\Vert
_{H_{p}}\rightarrow 0,\text{ when }n\rightarrow \infty .  \label{cnod1}
\end{equation}

b) For \ any $p\in \left( 0,1/2\right) $ there exists a martingale $f\in
H_{p}(G)$ \ \ for which
\begin{equation}
\omega _{H_{p}}\left( 1/2^{k},f\right) =O\left( 1/2^{k\left( 1/p-2\right)
}\right) ,\text{ \ as }k\rightarrow \infty .  \label{cond2}
\end{equation}%
and
\begin{equation*}
\left\Vert \sigma _{n}^{\kappa }\left( f\right) -f\right\Vert _{L_{p,\infty
}}\nrightarrow 0,\,\,\,\text{as\thinspace \thinspace \thinspace }%
n\rightarrow \infty .
\end{equation*}
\end{theorem}

\begin{theorem}
a) Let
\begin{equation}
\omega _{H_{1/2}}\left( 1/2^{k},f\right) =o\left( 1/k^{2}\right) ,\text{ \
as }k\rightarrow \infty .  \label{10A}
\end{equation}%
Then
\begin{equation*}
\left\Vert \sigma _{n}^{\kappa }\left( f\right) -f\right\Vert
_{H_{1/2}}\rightarrow 0,\text{ when }n\rightarrow \infty .
\end{equation*}

b) There exists a martingale $f\in H_{1/2}(G),$ for which
\begin{equation*}
\omega _{H_{1/2}}\left( 1/2^{k},f\right) =O\left( 1/k^{2}\right) ,\text{ \
as }k\rightarrow \infty .
\end{equation*}%
and
\begin{equation*}
\left\Vert \sigma _{n}\left( f\right) -f\right\Vert _{1/2}\nrightarrow
0\,\,\,\text{as\thinspace \thinspace \thinspace }n\rightarrow \infty .
\end{equation*}
\end{theorem}

\section{AUXILIARY PROPOSITIONS}

\begin{lemma}
\cite{we2} A martingale $f=\left( f_{n},\text{ }n\in \mathbb{N}\right) $ is
in $H_{p}\left( 0<p\leq 1\right) $ if and only if there exists a sequence $%
\left( a_{k},\text{ }k\in \mathbb{N}\right) $ of p-atoms and a sequence $%
\left( \mu _{k},\text{ }k\in \mathbb{N}\right) $ of real numbers such that
for every $n\in \mathbb{N}$
\end{lemma}

\begin{equation}
\qquad \sum_{k=0}^{\infty }\mu _{k}S_{M_{n}}a_{k}=f_{n},\text{ \ \ a.e.}
\label{2A}
\end{equation}%
\textit{and}%
\begin{equation*}
\text{\ \ \ }\sum_{k=0}^{\infty }\left\vert \mu _{k}\right\vert ^{p}<\infty .
\end{equation*}%
\textit{Moreover,} $\left\Vert f\right\Vert _{H_{p}}\backsim \inf \left(
\sum_{k=0}^{\infty }\left\vert \mu _{k}\right\vert ^{p}\right) ^{1/p},$
\textit{where the infimum is taken over all decomposition of }$f$\textit{\
of the form }(\ref{2A}).

\begin{lemma}
\cite{goginava1} Let $2<A\in \mathbb{N}$ and $%
q_{A}=2^{2A}+2^{2A-2}+...+2^{2}+2^{0}.$ Then
\begin{equation*}
q_{A-1}\left\vert K_{q_{A-1}}^{w}\left( x\right) \right\vert \geq
2^{2m+2s-3},
\end{equation*}%
for $x\in I_{2A}\left(
0,...,0,x_{2m}=1,0...,0,x_{2s}=1,x_{2s+1},...,x_{2A-1}\right) ,$ $%
m=0,1,...,A-3,$ $s=m+2,$ $m+3,...,A-1.$
\end{lemma}

\section{Proof of the Theorems}

\textbf{Proof of Theorem 1. }Let $0<p<1/2.$ Combining (\ref{c1}) and (\ref%
{5.1}) we conclude that
\begin{eqnarray}
&&\left\Vert \sigma _{n}^{\kappa }f\right\Vert
_{H_{p}}^{p}=\int_{G}\left\Vert \widetilde{\sigma _{n}^{\kappa }f^{\left(
t\right) }}\right\Vert _{p}^{p}dt  \label{5.3} \\
&=&\int_{G}\left\Vert \sigma _{n}^{\kappa }\widetilde{f^{\left( t\right) }}%
\right\Vert _{p}^{p}dt\leq c_{p}n^{1-2p}\int_{G}\left\Vert \widetilde{%
f^{\left( t\right) }}\right\Vert _{H_{p}}^{p}dt  \notag \\
&\sim &c_{p}n^{1-2p}\int_{G}\left\Vert f\right\Vert
_{H_{p}}^{p}dt=c_{p}n^{1-2p}\left\Vert f\right\Vert _{H_{p}}^{p}.  \notag
\end{eqnarray}%
Let $2^{m}<n\leq 2^{m+1}.$ Using (\ref{5.3}) we have that

\begin{eqnarray*}
&&\left\Vert \sigma _{n}^{\kappa }f-f\right\Vert _{H_{p}}^{p} \\
&\leq &\left\Vert \sigma _{n}^{\kappa }f-\sigma _{n}^{\kappa
}S_{2^{m}}f\right\Vert _{H_{p}}^{p}+\left\Vert \sigma _{n}^{\kappa
}S_{2^{m}}f-S_{2^{m}}f\right\Vert _{H_{p}}^{p}+\left\Vert
S_{2^{m}}f-f\right\Vert _{H_{p}}^{p} \\
&=&\left\Vert \sigma _{n}^{\kappa }\left( S_{2^{m}}f-f\right) \right\Vert
_{H_{p}}^{p}+\left\Vert S_{2^{m}}f-f\right\Vert _{H_{p}}^{p}+\left\Vert
\sigma _{n}^{\kappa }S_{2^{m}}f-S_{2^{m}}f\right\Vert _{H_{p}}^{p} \\
&\leq &c_{p}\left( n^{1-2p}+1\right) \omega _{H_{p}}^{p}\left(
1/2^{m},f\right) +\left\Vert \sigma _{n}^{\kappa
}S_{2^{m}}f-S_{2^{m}}f\right\Vert _{H_{p}}^{p}.
\end{eqnarray*}%
By simple calculation we get that%
\begin{eqnarray*}
&&\sigma _{n}^{\kappa }S_{2^{m}}f-S_{2^{m}}f=\frac{2^{m}}{n}\left(
S_{2^{m}}\sigma _{2^{m}}f-S_{2^{m}}f\right) \\
&=&\frac{2^{m}}{n}S_{2^{m}}\left( \sigma _{2^{m}}f-f\right) .
\end{eqnarray*}

Let $p>0.$ By applying (\ref{31}) we obtain that%
\begin{eqnarray}
&&\left\Vert \sigma _{n}^{\kappa }S_{2^{m}}f-S_{2^{m}}f\right\Vert
_{H_{p}}^{p}  \label{10} \\
&\leq &\frac{2^{pm}}{n^{p}}\left\Vert S_{2^{m}}\left( \sigma
_{2^{m}}f-f\right) \right\Vert _{H_{p}}^{p}  \notag \\
&\leq &\left\Vert \sigma _{2^{m}}f-f\right\Vert _{H_{p}}^{p}\rightarrow 0,%
\text{ when }k\rightarrow \infty \text{.}  \notag
\end{eqnarray}

It follows that if
\begin{equation*}
\omega _{H_{p}}\left( 1/2^{m},f\right) =o\left( 1/2^{m\left( 1/p-2\right)
}\right) ,\text{ as }m\rightarrow \infty .
\end{equation*}%
then
\begin{equation*}
\left\Vert \sigma _{n}^{\kappa }f-f\right\Vert _{H_{p}}\rightarrow 0,\text{
when }n\rightarrow \infty .
\end{equation*}%
Now, we are ready to prove second part of Theorem 1. We set
\begin{equation*}
a_{i}\left( x\right) =2^{i\left( 1/p-1\right) }\left( D_{2^{i+1}}\left(
x\right) -D_{2^{i}}\left( x\right) \right) ,
\end{equation*}%
and
\begin{equation*}
f^{\left( n\right) }\left( x\right) =\underset{i=0}{\overset{n}{\sum }}\frac{%
1}{2^{\left( 1/p-2\right) i}}a_{i}(x),\text{ }n\in \mathbb{N}.
\end{equation*}%
It is easy to show that $f=\left( f^{\left( 1\right) },f^{\left( 2\right)
},...,f^{\left( n\right) },...\right) \in H_{p}.$ Indeed, since

\begin{equation}
S_{2^{A}}a_{i}\left( x\right) =\left\{
\begin{array}{ll}
a_{i}\left( x\right) , & i<A, \\
0, & i\geq A,%
\end{array}%
\right.  \label{19}
\end{equation}

\begin{equation*}
\text{supp}(a_{i})=I_{i},\text{\ \ \ }\int_{I_{i}}a_{i}d\mu =0\text{\ }
\end{equation*}%
and

\begin{equation*}
\text{\ }\left\Vert a_{i}\right\Vert _{\infty }\leq 2^{\left( 1/p-1\right)
i}\cdot 2^{i}=2^{i/p}=\left( \text{supp }a_{i}\right) ^{-1/p},
\end{equation*}%
if we apply Lemma 1 we conclude that $f\in H_{p}.$

It is easy to show that
\begin{eqnarray}
&&f-S_{2^{n}}f  \label{20} \\
&=&\left( f^{\left( 1\right) }-S_{2^{n}}f^{\left( 1\right) },...,f^{\left(
n\right) }-S_{2^{n}}f^{\left( n\right) },...,f^{\left( n+k\right)
}-S_{2^{n}}f^{\left( n+k\right) }\right)  \notag \\
&=&\left( 0,...,0,f^{\left( n+1\right) }-f^{\left( n\right) },...,f^{\left(
n+k\right) }-f^{\left( n\right) },...\right)  \notag \\
&=&\left( 0,...,0,\underset{i=n}{\overset{n+k}{\sum }}\frac{a_{i}(x)}{%
2^{i\left( 1/p-2\right) }},...\right) ,\text{ \ }k\in \mathbb{N}_{+}.  \notag
\end{eqnarray}%
is a martingale. By using (\ref{20}) we get that%
\begin{equation*}
\omega _{H_{p}}(1/2^{n},f)
\end{equation*}%
\begin{equation*}
\leq \sum\limits_{i=n+1}^{\infty }1/2^{\left( 1/p-2\right) i}\leq
c/2^{\left( 1/p-2\right) n}=O\left( 1/2^{\left( 1/p-2\right) n}\right) .
\end{equation*}%
It is easy to show that

\begin{equation}
\widehat{f}(j)=2^{i},\text{ \ \ if \thinspace \thinspace }j\in \left\{
2^{i},...,2^{i+1}-1\right\} ,\text{ }i=0,1,...  \label{29}
\end{equation}

By applying (\ref{31}) and (\ref{29}) we have that
\begin{eqnarray*}
&&\limsup\limits_{n\rightarrow \infty }\Vert \sigma _{2^{n}+1}^{\kappa
}(f)-f\Vert _{L_{p,\infty }} \\
&=&\limsup\limits_{n\rightarrow \infty }\Vert \frac{2^{n}\sigma
_{2^{n}}^{\kappa }(f)}{2^{n}+1}+\frac{S_{2^{n}}(f)}{2^{n}+1}+\frac{%
2^{n}\kappa _{2^{n}}}{2^{n}+1}-\frac{2^{n}f}{2^{n}+1}-\frac{f}{2^{n}+1}\Vert
_{L_{p,\infty }} \\
&\geq &\limsup\limits_{n\rightarrow \infty }\frac{2^{n}}{2^{n}+1}\Vert
\kappa _{2^{n}}\Vert _{L_{p_{,\infty }}} \\
&&-\limsup\limits_{n\rightarrow \infty }\frac{2^{n}}{2^{n}+1}\Vert \sigma
_{2^{n}}^{\kappa }(f)-f\Vert _{L_{p,\infty }} \\
&&-\limsup\limits_{n\rightarrow \infty }\frac{1}{2^{n}+1}\Vert
S_{2^{n}}(f)-f\Vert _{L_{p,\infty }} \\
&>&\limsup\limits_{n\rightarrow \infty }\frac{2^{n}}{2^{n}+1}-o(1)>0.
\end{eqnarray*}

Theorem 1 is proved.

\textbf{Proof of Theorem 2. }Analogously of first part of Theorem 1 we can
show that if modulus of continuity of the martingale $f\in H_{1/2}$
satisfies the condition (\ref{10A}), then $\left\Vert \sigma _{n}^{\kappa
}f-f\right\Vert _{H_{1/2}}\rightarrow \infty ,$ when $n\rightarrow \infty .$

Now, we are ready to prove the second part of theorem 1. We set
\begin{equation*}
a_{i}(x)=2^{2^{i}}\left( D_{2^{2^{i}+1}}(x)-D_{2^{2^{i}}}(x)\right)
\end{equation*}%
and
\begin{equation*}
f^{\left( n\right) }\left( x\right) =\overset{n}{\sum_{i=1}}\frac{a_{i}(x)}{%
2^{2i}}.
\end{equation*}

It is easy to show that $\ f=\left( f^{\left( 1\right) },f^{\left( 2\right)
},...,f^{\left( n\right) },...\right) \in H_{1/2}.$ Indeed, since

\begin{equation*}
S_{2^{A}}a_{i}\left( x\right) =\left\{
\begin{array}{ll}
a_{i}\left( x\right) , & 2^{i}<A, \\
0, & 2^{i}\geq A%
\end{array}%
\right.
\end{equation*}

\begin{equation*}
\text{supp}(a_{i})=I_{2^{i}},\text{ \ \ }\int_{I_{2^{i}}}a_{i}d\mu =0,\text{
\ \ \ }\left\Vert a_{i}\right\Vert _{\infty }\leq \left( 2^{2^{i}}\right)
^{2}=\left( \mu (\text{supp }a_{i})\right) ^{2},
\end{equation*}%
if we apply Lemma 1 we conclude that $f\in H_{1/2}.$

By simple calculation we get that%
\begin{equation*}
\omega _{H_{1/2}}\left( 1/2^{n},f\right) \leq \overset{\infty }{\sum_{i=%
\left[ \log _{2}^{n}\right] }}\frac{1}{2^{2i}}=O\left( 1/n^{2}\right) ,
\end{equation*}%
where $\left[ \log _{2}^{n}\right] $ denotes the integer part of $\log
_{2}^{n}.$ Hence%
\begin{eqnarray}
&&\Vert \sigma _{q_{2^{k-1}}}^{\kappa }(f)-f\Vert _{1/2}  \label{nn} \\
&=&\Vert \frac{2^{2^{k}}\sigma _{2^{2^{k}}}^{\kappa }(f)}{q_{_{2^{k-1}}}}+%
\frac{1}{q_{_{2^{k-1}}}}\sum_{j=2^{2^{k}}+1}^{q_{2^{k-1}}}S_{j}^{\kappa }f-%
\frac{2^{2^{k}}f}{q_{_{2^{k-1}}}}-\frac{q_{_{2^{k-1}-1}}f}{q_{_{2^{k-1}}}}%
\Vert _{1/2}  \notag
\end{eqnarray}

It is easy to show that

\begin{equation}
\widehat{f}(j)=\left\{
\begin{array}{ll}
\frac{2^{2^{i}}}{2^{2i}}, & \text{if \thinspace \thinspace }j\in \left\{
2^{2^{i}},...,2^{2^{i}+1}-1\right\} ,\text{ }i=0,1,... \\
0, & \text{\thinspace if \thinspace \thinspace \thinspace }j\notin
\bigcup\limits_{n=0}^{\infty }\left\{ 2^{2^{n}},...,2^{2^{n}+1}-1\right\} .%
\text{ }%
\end{array}%
\right.  \label{35}
\end{equation}

Let $2^{2^{i}}<j\leq q_{_{2^{i-1}}}.$ Since
\begin{equation*}
D_{s+2^{2^{i}}}^{\kappa }=D_{2^{2^{i}}}+r_{2^{^{i}}}D_{s}^{w}\left( \tau
_{2^{i}}\right) ,\text{ \qquad when \thinspace \thinspace }s<2^{2^{i}}
\end{equation*}%
by applying (\ref{35}) we obtain that
\begin{eqnarray*}
S_{j}^{\kappa }f &=&S_{2^{2^{i}}}f+\sum_{v=2^{2^{i}}}^{j-1}\widehat{f}%
(v)\kappa _{v}=S_{2^{2^{i}}}f+\frac{2^{2^{i}}}{2^{2i}}\sum_{v=2^{i}}^{j-1}%
\kappa _{v} \\
&=&S_{2^{2^{i}}}f+\frac{2^{2^{i}}}{2^{2i}}\left( D_{_{j}}^{\kappa
}-D_{2^{2^{i}}}\right) =S_{2^{2^{i}}}f+\frac{%
2^{2^{i}}r_{2^{^{i}}}D_{_{j-2^{2^{i}}}}^{w}\left( \tau _{2^{i}}\right) }{%
2^{2i}}.
\end{eqnarray*}

Hence
\begin{eqnarray*}
&&\sum_{j=2^{2^{i}}+1}^{q_{2^{i-1}}}S_{j}^{\kappa }f=\frac{%
q_{2^{i-1}-1}S_{2^{2^{i}}}f}{q_{2^{i-1}}}+\frac{2^{2^{i}}r_{2^{^{i}}}}{%
q_{2^{i}}2^{2i}}\sum_{j=2^{2^{i}}+1}^{q_{2^{i-1}}}D_{_{j-2^{2^{k}}}}^{w}%
\left( \tau _{2^{i}}\right) \\
&=&\frac{q_{2^{i-1}-1}S_{2^{2^{i}}}f}{q_{2^{i-1}}}+\frac{%
2^{2^{i}}r_{2^{^{i}}}}{q_{2^{i-1}}2^{2i}}%
\sum_{j=1}^{q_{2^{i-1}-1}}D_{_{j}}^{w}\left( \tau _{2^{i}}\right) \\
&=&\frac{q_{2^{i-1}-1}S_{2^{2^{i}}}f}{q_{2^{i-1}}}+\frac{%
2^{2^{i}}r_{2^{^{i}}}q_{2^{i-1}-1}K_{q_{2^{i-1}-1}}^{w}\left( \tau
_{2^{i}}\right) }{q_{2^{i-1}}2^{2i}}.
\end{eqnarray*}

By using (\ref{nn}) we get that%
\begin{eqnarray}
&&\Vert \sigma _{q_{2^{i-1}}}^{\kappa }(f)-f\Vert _{1/2}^{1/2}  \label{a11}
\\
&\geq &\frac{c}{2^{i}}\Vert q_{2^{i-1}-1}K_{q_{2^{i-1}-1}}^{w}\left( \tau
_{2^{i}}\right) \Vert _{1/2}^{1/2}-\left( \frac{2^{2^{i}}}{q_{2^{i-1}}}%
\right) ^{1/2}\Vert \sigma _{_{2^{2^{i}}}}f-f\Vert _{1/2}^{1/2}  \notag \\
&&-\left( \frac{q_{2^{i-1}-1}}{q_{2^{i-1}}}\right) ^{1/2}\Vert
S_{2^{2^{i}}}f-f\Vert _{1/2}^{1/2}=IV-V-VI.  \notag
\end{eqnarray}%
Let $x\in J_{2^{i}}^{2^{i}-2s-1,2^{i}-2\eta -1},$ where $\eta
=1,...,2^{i}-3,\,s=\eta +3,...,2^{i}-1.$ By applying Lemma 2 we have that

\begin{equation*}
q_{2^{i-1}-1}\left\vert K_{q_{2^{i-1}-1}}^{w}\left( \tau _{2^{i}}\left(
x\right) \right) \right\vert \geq 2^{2\eta +2s-3}.
\end{equation*}

We can write

\begin{eqnarray*}
&&\int_{G}\left\vert q_{2^{i-1}-1}K_{q_{2^{i-1}-1}}^{w}\left( \tau
_{2^{i}}\left( x\right) \right) \right\vert ^{1/2}d\mu \left( x\right) \\
&\geq &c\sum_{\eta =1}^{2^{i}-3}\sum_{s=\eta
+2}^{2^{i}-1}\sum_{x_{2s+1}=0}^{1}...\sum_{x_{_{2^{i}-1}}=0}^{1}%
\int_{J_{2^{i}}^{2^{i}-2s-1,2^{i}-2\eta -1}}\left\vert
q_{2^{i-1}-1}K_{q_{2^{i-1}-1}}^{w}\left( \tau _{2^{i}}\left( x\right)
\right) \right\vert ^{1/2}d\mu \left( x\right) \\
&\geq &c\sum_{\eta =1}^{2^{i}-3}\sum_{s=\eta +2}^{2^{i}-1}\frac{1}{2^{2s}}%
\sqrt{2^{2\eta +2s}}\geq \sum_{\eta =1}^{2^{i}-3}\sum_{s=\eta +2}^{2^{i}-1}%
\sqrt{\frac{2^{2\eta }}{2^{2s}}}\geq c\sum_{\eta =1}^{2^{i}-3}\geq c2^{i}.
\end{eqnarray*}

By using (\ref{a11}) we have that%
\begin{equation*}
\limsup\limits_{i\rightarrow \infty }\Vert \sigma _{q_{2^{i-1}}}^{\kappa
}(f)-f\Vert _{1/2}\geq c>0.
\end{equation*}

Theorem 2 is proved.

\textbf{Acknowledgment: }The author would like to thank the referees for
their helpful suggestions.

\end{document}